\documentclass[a4paper,11pt]{article}

\usepackage{fullpage}
\usepackage{amsmath,amsthm,amssymb, amsfonts, mathrsfs}
\usepackage{mathtools}
\usepackage[shortlabels]{enumitem}
\usepackage[hidelinks]{hyperref}
\usepackage[nameinlink,capitalise,noabbrev]{cleveref}
\usepackage{graphicx, float, tikz, subcaption} 
\usepackage{multirow}

\theoremstyle{definition}
\newtheorem{definition}{Definition}[section]

\theoremstyle{plain}

\newtheorem{theorem}[definition]{Theorem}

\def \cl {\colon}
\def \ce {\coloneqq}

\renewcommand{\le}{\leqslant}
\renewcommand{\ge}{\geqslant}

\makeatletter

\newcommand*{\rom}[1]{\expandafter\@slowromancap\romannumeral #1@}

\def \mA {\mathcal{A}}

\def \mF{\mathcal{F}}
\def \mG {\mathcal{G}}

\newcommand{\floor}[1]{\left \lfloor #1 \right \rfloor}

\DeclareMathOperator{\VC}{VC}
\renewcommand {\vec}[1]{\overrightarrow{#1}}

\title{VC-dimensions Between Partially Ordered Sets and Totally Ordered Sets}

\author{
	Boyan Duan\thanks{School of Computer Science, ETH Z\"urich, Z\"urich 8092, Switzerland. \texttt{boduan@student.ethz.ch}. }
	\and
	Minghui Ouyang\thanks{School of Mathematical Sciences, Peking University, Beijing 100871, China. \texttt{ouyangminghui1998@gmail.com}. }
	\and
	Zheng Wang\thanks{School of Mathematical Sciences, Peking University, Beijing 100871, China. \texttt{wangzheng0401@hotmail.com}. }
}

\date{\vspace{-5ex}}

\begin{document}
	\maketitle
	\begin{abstract}
		We say that two partial orders on $[n]$ are compatible if there exists a partial order that refines both of them. This compatibility relation induces a natural set system structure between the collection $\mathcal{F}$ of all partial orders and the collection $\mathcal{G}$ of all total orders on $[n]$, where each order is associated with the set of orders compatible with it. 
		
		In this note, we determine the VC-dimension of $\mathcal{F}$ with respect to $\mathcal{G}$, proving that $\operatorname{VC}_{\mathcal{G}}(\mathcal{F}) = \lfloor\frac{n^2}{4}\rfloor$ for $n \ge 4$. We also establish bounds on the dual VC-dimension, showing that $2(n-3) \le \operatorname{VC}_{\mathcal{F}}(\mathcal{G}) \le n \log_2 n$ for all $n \ge 1$. 
		
	\end{abstract}
	
\section{Introduction}
	For a family $\mF$ of subsets of a set $X$, a subset $S \subseteq X$ is said to be \emph{shattered} by $\mF$ if, for every $A \subseteq S$ there exists $B \in \mF$ such that $B \cap S = A$. The \emph{VC-dimension} of $\mF$ is the largest cardinality of a subset of $X$ that is shattered by $\mF$. We denote the VC-dimension of $\mF$ with respect to $X$ by $\VC_{X}(\mF)$. Since its introduction by Vapnik and Chervonenkis~\cite{VC71} in the context of statistical learning theory, VC-dimension has played a central role in various areas of mathematics and computer science. 

	If $\mF$ shatters $S$, then each subset of $S$ corresponds to a distinct element of $\mF$. In particular, the VC-dimension of a set system $\mF$ is always at most $\log_2 |\mF|$, regardless of the ground set. 

	In this paper we fix the ground set $[n]$ and consider the collections of partial and total orders on $[n]$, studied under the following notion of \emph{compatibility}. 

	\begin{definition} \label{def:compatibility}
		Given $n \in \mathbb{N}$, let $\mF$ (resp. $\mG$) denote the set of all partial (resp. total) orders on $[n]$. Clearly, $\mG \subseteq \mF$. We say that two partial orders $<_1$ and $<_2$ on $[n]$ are \emph{compatible} if there exists a partial order that refines both $<_1$ and $<_2$. 
	\end{definition}

	Equivalently, one may view a partial order as a directed graph on vertex set $[n]$, where an edge $a \to b$ is drawn whenever $a < b$ in the order. Then $<_1$ and $<_2$ are compatible if and only if the directed graph formed by $<_1 \cup <_2$ is acyclic. In the special case where $<_1$ is a total order, this is equivalent to requiring that $<_1$ be a linear extension of $<_2$. 

	Under this compatibility relation, $\mF$ and $\mG$ naturally induce set systems on each other. 

	\begin{definition} \label{def:set_families}
		For each element $A \in \mF$ (resp. $A \in \mG$), we define the \emph{associated} subset of $\mG$ (resp. $\mF$) to be the collection of all $B \in \mG$ (resp. $B \in \mF$) that are compatible with $A$. Accordingly, we define $\VC_{\mG}(\mF)$ as the VC-dimension of $\mF$ when viewed as the associated set system over $\mG$. The quantities $\VC_{\mF}(\mF)$, $\VC_{\mF}(\mG)$ and $\VC_{\mG}(\mG)$ are defined analogously. 
	\end{definition}
	
	Kleitman and Rothschild~\cite{KR70} used a sophisticated counting argument\footnote{More precisely, they classified the set of partial orders $\mF$ into several classes and obtained a sharper bound $\log_2 |\mF| \le n^2/4 + O(n^{3/2} \log n)$. In subsequent work~\cite{KR75, BPS96}, the asymptotic value of $|\mF|$ was determined. However, these results are still insufficient for determining the VC-dimension. } to show that the size of $\mF$ is approximately $2^{(1+o(1)) \frac{n^2}{4}}$. Consequently, we obtain the upper bound 
	\[ \VC_{\mG}(\mF) \,\le\, \VC_{\mF}(\mF) \,\le\, \log_2 2^{(1+o(1))\,\tfrac{n^2}{4}} \,=\, (1+o(1))\,\tfrac{n^2}{4}. \] 
	In this note, we show that the $(1+o(1))$ factor can in fact be removed. 
	
	\begin{theorem} \label{thm:vc_dim_partial_to_total}
		For $n \ge 1$, we have
		\[ \VC_{\mG}(\mF) \,=\,
		\begin{cases}
			3,& n = 3, \\
			\floor{\tfrac{n^2}{4}} ,& n \neq 3. 
		\end{cases}
		\]
	\end{theorem}
	
	We believe that $\VC_{\mF}(\mF)$ is also equal to $\floor{\tfrac{n^2}{4}}$, although our proof does not extend directly to this case. 
	
	For the VC-dimension of $\mG$, it is easy to see that $\VC_{\mG}(\mG) = 1$ when $\mG$ is viewed as a set system on itself, since the only total order compatible with a given total order $<$ is $<$ itself. Thus, the only remaining meaningful question is to determine the value of $\VC_{\mF}(\mG)$. 
	
	We establish the following lower bound by explicitly constructing a shattered set. We believe, however, that the correct asymptotic value of $\VC_{\mF}(\mG)$ is $\Theta(n \log n)$. 
	
	\begin{theorem} \label{thm:vc_dim_total_to_partial} 
		For $n \ge 1$, we have
			\[ 2(n-3) \,\le\, \VC_{\mF}(\mG) \,\le\, n \log_2 n. \]
	\end{theorem}
	
\section{Proofs}

	\begin{proof}[\underline{Proof of \Cref{thm:vc_dim_partial_to_total}}]
		For $n = 1,2$, it is easy to check that $\VC_\mG(\mF) = n-1 = \lfloor \tfrac{n^2}{4} \rfloor$. For $n = 3$, consider the set $S = \{123, 231, 312\}$, where the string ``$abc$'' denotes the total order $a < b < c$ on $\{a, b, c\}$. We claim that $S$ is shattered, yielding that the VC-dimension is at least $3$. Indeed, for each subset of $S$, one can find a poset whose set of compatible total orders intersects $S$ exactly in that subset. 
		\vspace{-1ex}\begin{itemize}[itemsep=0pt, parsep=0pt]
			\item The empty poset $\varnothing$ is compatible with all elements of $S$. 
			\item The poset $23$ (i.e., $2 < 3$ and $1$ is incomparable with $2$ and $3$) is compatible with $123$ and $231$, but not with $312$. 
			\item The poset $123$ is compatible only with $123$, and incompatible with $231$ and $312$. 
			\item The poset $321$ is incompatible with all elements of $S$. 
		\end{itemize}
	
		Thus $S$ is shattered, and hence $\VC_\mG(\mF) \ge 3$. To see that the VC-dimension is less than $4$, note that the only poset compatible with both $abc$ and $cba$ is the empty poset $\varnothing$, which is in fact compatible with every total order. Therefore, if a shattered set contained a pair of reverse orders, it could consist of at most those two elements. Consequently, no shattered set of size $4$ exists, and we conclude that $\VC_\mG(\mF) = 3$ when $n=3$. 
	
		Now assume $n \ge 4$. We represent partial orders as directed graphs, where an edge $a \to b$ corresponds to the relation $a < b$ in the partial order. Suppose $S \subseteq \mG$ is a shattered set, we aim to show that $|S| \le \floor{\tfrac{n^2}{4}}$. Since $n \le \floor{\tfrac{n^2}{4}}$ for every $n \ge 4$, we may assume $|S| \ge n+1$ for contradiction, as otherwise the bound is already established. 
	
		By the definition of shattering, for each $A \in S$, there exists a partial order $G_A$ on $[n]$ that is incompatible with $A$ but compatible with all orders in $S \setminus \{ A \}$. Since $G_A$ is incompatible with $A$ and $A$ is a total order, there exists an edge $e_A \in G_A$ that contradicts $A$. Define $G \ce \{e_A \cl A \in S\}$ to be the directed graph formed by the set of edges $e_A$. We establish the following two claims: 
	
		\vspace{1ex}
		\medskip
		\textbf{Claim 1: $G$ is acyclic.} 
		\begin{proof}
			Suppose, for contradiction, that $G$ contains a directed cycle $C = e_{A_1} e_{A_2} \cdots e_{A_k}$. Without loss of generality, assume that $C$ is a simple cycle, so $k \le n$. Since $|S| \ge n+1$, there exists some $B \in S \setminus \{A_1, \cdots, A_k\}$. 
		
			By the definition of $e_{A_i}$, each edge $e_{A_i}$ is contained in $G_{A_i}$, and is therefore compatible with all orders in $S$ except $A_i$. In particular, $e_{A_i}$ is compatible with $B$. 
		
			Because $B$ is a total order, it determines the relative order of every pair of elements. Thus, each edge $e_{A_i}$ must lie in $B$. Therefore, the cycle $C$ is entirely contained in $B$, contradicting the fact that a total order cannot contain a directed cycle. 
		\end{proof}
	
		\medskip
		\textbf{Claim 2: $G$ does not contain a directed path of length at least $2$ from $x$ to $y$ if the edge $\vec{xy} \in G$.}  
		\begin{proof}
			Assume for contradiction that there exists a directed path $P = e_{A_1} e_{A_2} \cdots e_{A_k}$ from $x$ to $y$ with $k \ge 2$, and that the edge $\vec{xy}$ also appears in $G$ as $e_B = \vec{xy}$ for some $B \in S$. 
		
			By the same argument as in Claim 1, each edge $e_{A_i}$ must belong to $B$. Since $B$ is a total order, it follows that $\vec{xy} \in B$. However, this contradicts the fact that $\vec{xy} = e_B$ was chosen to contradict $B$. Therefore, no such path can exist. 
		\end{proof}
	
		From the above two claims, we conclude that $G$ contains no directed cycle and no directed path of length $\ge 2$ between any pair $(x,y)$ such that $\vec{xy} \in G$. In particular, the underlying undirected graph of $G$ is triangle-free. Therefore, $|S| = |E(G)| \le \floor{\tfrac{n^2}{4}}$, as desired. 
	
		Next, we show that there exists a shattered set $S \subseteq \mG$ of size $\floor{\tfrac{n^2}{4}}$. Consider all pairs $(i,j)$ such that $1 \le i \le \floor{\tfrac{n}{2}}$ and $\floor{\tfrac{n}{2}}+1 \le j \le n$. For each such pair, define the directed edge $e_{i,j} = \vec{ij}$ and its reverse $\overline{e_{i,j}} = \overleftarrow{ij}$. 
	
		For each pair $(i,j)$, let $A_{i,j}$ be a topological ordering of the acyclic graph consisting of $\{\overline{e_{i,j}}\} \cup \{e_{i',j'} \cl (i',j') \neq (i,j)\}$ (that is, $A_{i,j}$ is a linear extension of the partial order defined by this directed acyclic graph). We claim that the set 
			\[ S \,\ce\, \biggl\{A_{i,j} \cl 1 \le i \le \floor{\tfrac{n}{2}}, \floor{\frac{n}{2}}+1 \le j \le n \biggr\} \] 
		is shattered. 
	
		To see this, consider any subset $\mA \subseteq S$. Define a partial order $G = \{ e_{i,j} \cl A_{i,j} \in \mA\}$. By construction, $G$ is compatible with elements of $S \setminus \mA$ but incompatible with elements of $\mA$, proving that $S$ is a shattered set. 
	\end{proof}
	
	\begin{proof}[\underline{Proof of \Cref{thm:vc_dim_total_to_partial}}]
		Since $\mG$ consists of all total orders on $[n]$, we have $|\mG| = n!$. Therefore,
			\[ \VC_{\mF}(\mG) \,\le\, \log_2 |\mG| \,\le\, n \log_2 n. \]
		
		To establish a lower bound for $\VC_{\mF}(\mG)$, we present two distinct constructions of shattered sets, yielding
			\[ \VC_{\mF}(\mG) \,\ge\, 3 \Bigl(\left\lfloor \tfrac{n}{2} \right\rfloor - 1 \Bigr) \quad\text{and}\quad \VC_{\mF}(\mG) \,\ge\, 2(n-3), \]
		respectively. Since these two constructions rely on rather different structures, we describe both.

		Each construction produces a family of directed graphs, denoted $\{P_i\}$ and $\{Q_i\}$, whose transitive closures (viewed as partial orders) form a shattered set. 
		
		For any subset $S \subseteq \{P_i\}$ (or $S \subseteq \{Q_i\}$), there exists a way to select one edge from each graph in $S$ and delete all other edges, such that if we reverse the direction of the selected edges while leaving all other graphs from $\{P_i\} \setminus S$ (resp. $\{Q_i\} \setminus S$) unchanged, the resulting union of edges is acyclic: 
			\[ \forall S \subseteq \{P_i\},\, \exists\, \{e_i \in E(P_i) \colon P_i \in S\} \quad \text{such that} \quad \Bigl( \bigcup_{P_i \in S} \{\overline{e_i}\} \Bigr) \,\cup\, \Bigl( \bigcup_{P_i \notin S} E(P_i) \Bigr) \,\, \text{is acyclic}. \qquad (\boldsymbol{\ast}) \]
		Here $\overline{e_i}$ denotes the edge $e_i$ with its direction reversed.
		
		We refer to this as property ($\boldsymbol{\ast}$). This condition is precisely equivalent to saying that the transitive closures of $\{P_i\}$ (resp. $\{Q_i\}$) form a shattered set: for every $S$, any topological ordering of the resulting acyclic graph yields a total order that is compatible with all graphs outside $S$ and incompatible with all graphs inside $S$. 
		
		\begin{figure}[ht!]
			\centering
			\begin{subfigure}{.5\textwidth}
				\centering
				\includegraphics[width=.9\linewidth]{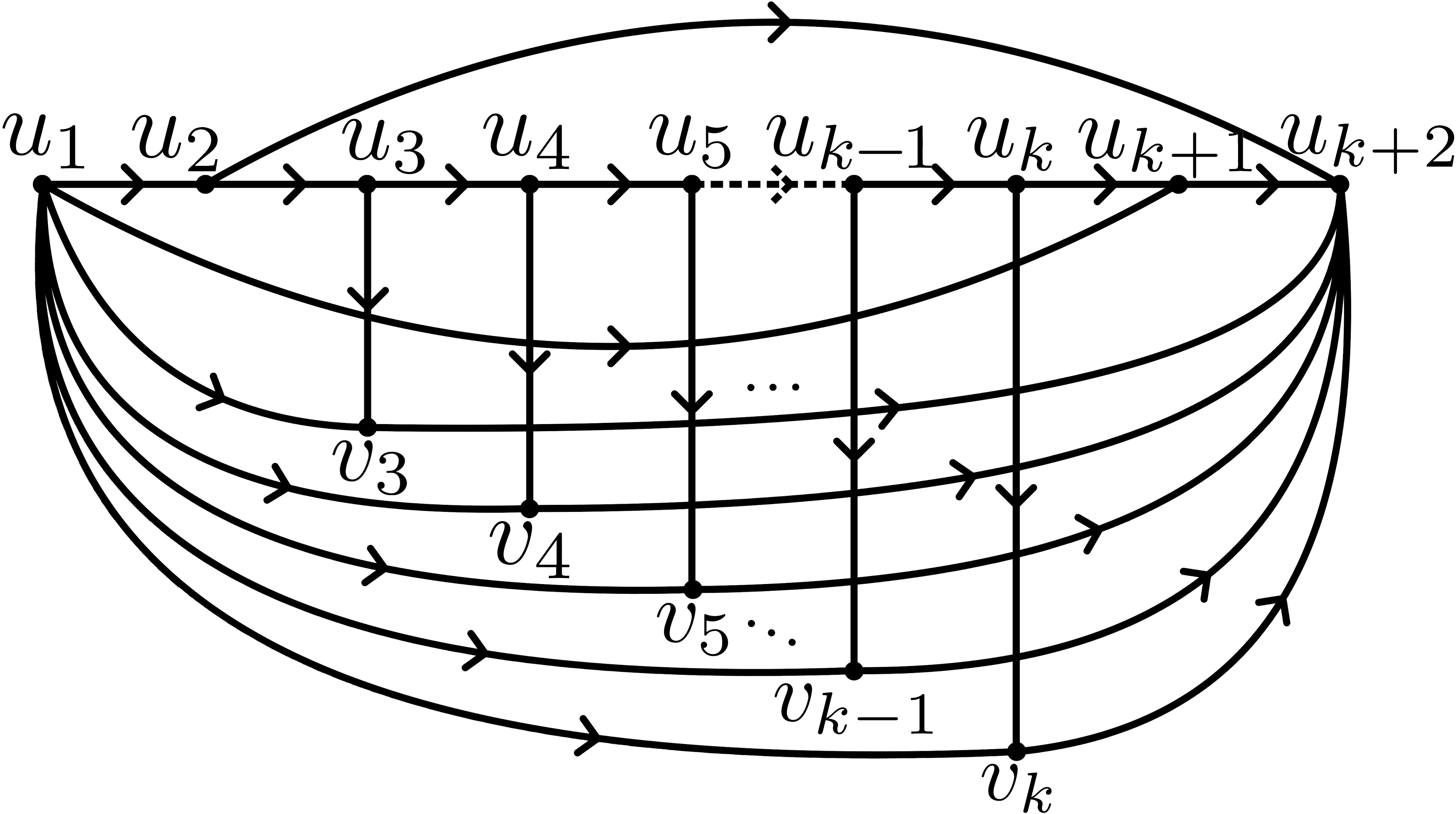}
				\caption{The graph $P$}
			\end{subfigure}%
			\begin{subfigure}{.5\textwidth}
				\centering
				\includegraphics[width=.8\linewidth]{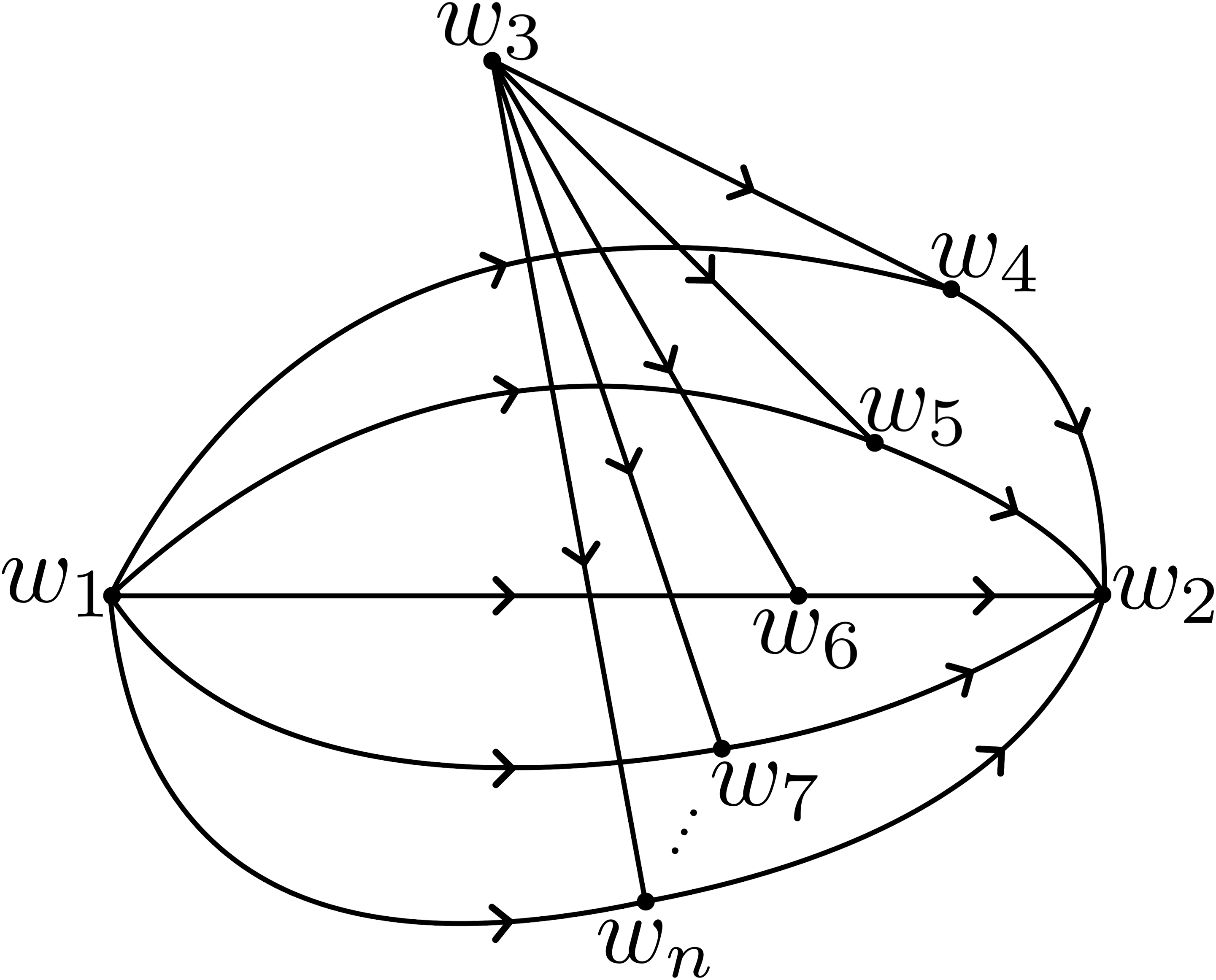}
				\caption{The graph $Q$}
			\end{subfigure}
		\end{figure}
		
		We now proceed to describe the families $\{P_i\}$ and $\{Q_i\}$. Let $k \ce \floor{\frac{n}{2}}$. We define
		\begin{gather*}
			P_1 \ce \bigl\{ \vec{u_1u_2}, \vec{u_2u_{k+2}} \bigr\},\ P_2 \ce \bigl\{ \vec{u_2u_3} \bigr\}, \cdots,\ P_k \ce \bigl\{ \vec{u_ku_{k+1}} \bigr\},\ P_{k+1} \ce \bigl\{ \vec{u_1u_{k+1}}, \vec{u_{k+1}u_{k+2}} \bigr\}, \\
			P_{k+2} \ce \bigl\{ \vec{u_1v_3}, \vec{v_3u_{k+2}} \bigr\}, \cdots,\ P_{2k-1} \ce \bigl\{ \vec{u_1v_k}, \vec{v_ku_{k+2}} \bigr\}, \\
			P_{2k} \ce \bigl\{ \vec{u_3v_3} \bigr\}, \cdots,\ P_{3k-3} \ce \bigl\{ \vec{u_kv_k} \bigr\}, 
		\end{gather*} 
		and
		\begin{gather*}
			Q_1 \ce \bigl\{ \vec{w_1w_4}, \vec{w_4w_2} \bigr\},\ Q_2 \ce \bigl\{ \vec{w_1w_5}, \vec{w_5w_2} \bigr\}, \cdots,\ Q_{n-3} \ce \bigl\{ \vec{w_1w_n}, \vec{w_nw_2} \bigr\}, \\
			Q_{n-2} \ce \bigl\{ \vec{w_3w_4} \bigr\}, \cdots,\ Q_{2n-6} \ce \bigl\{ \vec{w_3w_n} \bigr\}. 
		\end{gather*}
		
		We now verify \textbf{($\boldsymbol{\ast}$)} for the families $\{P_i\}$ and $\{Q_i\}$ respectively. It is straightforward to check that both $P \ce \bigcup_{i = 1}^{3k-3} P_i$ and $Q \ce \bigcup_{i = 1}^{2n-6} Q_i$ are acyclic directed graphs on $n$ vertices (with an additional isolated vertex added to $P$ when $n$ is odd), and that they decompose into $3\bigl(\lfloor n/2 \rfloor - 1\bigr)$ and $2(n-3)$ subgraphs, respectively. 
		
		\vspace{1ex}
		\underline{\textbf{$\boldsymbol{\{P_i\}}$ satisfies ($\boldsymbol{\ast}$):}} \hspace{.5em} 
		Suppose $S$ is an arbitrary subset of $\{P_i\}$. The graphs $P_2, \cdots, P_k$ and $P_{2k}, \cdots, P_{3k-3}$ each contain only a single edge, so their orientations after reversal are uniquely determined. If at least one of $P_2, \cdots, P_k$ belongs to $S$, then the directed path $u_2 \!\to\! u_3 \!\to\! \cdots \!\to\! u_{k+1}$ would be destroyed. In this case: 
		\vspace{-1ex}\begin{itemize}[itemsep=0pt, parsep=0pt]
			\item We reverse $\vec{u_2u_{k+2}}$ and delete $\vec{u_1u_2}$ if $P_1 \in S$, 
			\item and reverse $\vec{u_1u_{k+1}}$ and delete $\vec{u_{k+1}u_{k+2}}$ if $P_{k+1} \in S$. 
		\end{itemize}\vspace{-1ex}
		
		Otherwise:
		\vspace{-1ex}\begin{itemize}[itemsep=0pt, parsep=0pt]
			\item We reverse $\vec{u_1u_2}$ and delete $\vec{u_2u_{k+2}}$ if $P_1 \in S$, 
			\item and reverse $\vec{u_{k+1}u_{k+2}}$ and delete $\vec{u_1u_{k+1}}$ if $P_{k+1} \in S$. 
		\end{itemize}\vspace{-.5ex}
		
		From the figure of the graph $P$, the only two potential cycles within the vertex set $\{u_1, \cdots, u_{k+2}\}$, regardless of orientation, are $u_1 u_2 \cdots u_{k+1} u_1$ and $u_2 u_3 \cdots u_{k+2} u_2$. Under the above strategy, one can verify that neither of these cycles occurs in the resulting graph (as described in ($\boldsymbol{\ast}$)). 
		
		Therefore, any potential directed cycle in the resulting graph must involve vertices from $\{v_3, \cdots, v_k\}$. We then reverse the graphs in $S \cap \{P_{k+2}, \cdots, P_{2k-1}\}$ as follows: for each $P_i \in S \cap \{P_{k+2}, \dots, P_{2k-1}\}$,
		\vspace{-1ex}
		\begin{itemize}[itemsep=0pt, parsep=0pt]
			\item if $P_{i+k-2} \in S$, reverse $\vec{u_1 v_{i-k+1}}$ and delete $\vec{v_{i-k+1} u_{k+2}}$; 
			\item otherwise, reverse $\vec{v_{i-k+1} u_{k+2}}$ and delete $\vec{u_1 v_{i-k+1}}$.
		\end{itemize}\vspace{-.5ex} 
		
		Under this strategy, the only directed paths passing through $\{v_3, \cdots, v_k\}$ are of the form $u_1 \!\to\! v_i \!\to\! u_{k+2}$. Moreover, there is no directed path from $u_{k+2}$ back to $u_1$. It follows that these paths cannot form a directed cycle. 
		
		\vspace{1ex}
		\underline{\textbf{$\boldsymbol{\{Q_i\}}$ satisfies ($\boldsymbol{\ast}$):}} \hspace{.5em} 
		Suppose $S$ is an arbitrary subset of $\{Q_i\}$. Define $S_1 \ce S \cap \{Q_i \cl Q_{i+(n-3)} \in S\}$, $S_2 \ce S \cap \{Q_i \cl Q_{i+(n-3)} \notin S\}$ and $S_3 \ce S \cap \{Q_{n-2}, \cdots, Q_{2n-6}\}$. Then $S = S_1 \sqcup S_2 \sqcup S_3$ is a decomposition of $S$. We manipulate subgraphs in $S$ according to the following strategy: 
		\vspace{-1ex}\begin{itemize}[itemsep=0pt, parsep=0pt]
			\item For each $Q_i \in S_1$, we reverse $\vec{w_1w_{i+3}}$ and delete $\vec{w_{i+3}w_2}$;
			\item For each $Q_i \in S_2$, we reverse $\vec{w_{i+3}w_2}$ and delete $\vec{w_1w_{i+3}}$; 
			\item For each $Q_i \in S_3$, we reverse the only edge $\vec{w_3w_{i-(n-6)}}$. 
		\end{itemize}\vspace{-.5ex}
		
		One can verify that there is no directed path from $w_2$ to $w_1$ after the reversal. Therefore, any directed cycle in the resulting graph cannot simultaneously pass through $w_1$ and $w_2$. Hence, if a cycle exists, it must pass through both $w_1$ and $w_3$, or both $w_2$ and $w_3$. 
		
		However, under the above reversal strategy, there is no directed path from $w_3$ to $w_1$, nor from $w_2$ to $w_3$. It follows that the resulting graph contains no directed cycle, completing the proof. 
	\end{proof}
	
	\bibliographystyle{abbrv}
	\bibliography{reference}
	
\end{document}